\documentclass[12pt]{article}
\usepackage{amssymb}

\begin {document}

\title{Mathematical formulae of the all-interval tone-semitone series and the quart modes}

\author{Grigori G. Amosov \footnote {E-mail:gramos@mail.ru}\, and Sergei V. Golubkov \footnote {E-mail:sgolubk@mail.ru}\\
\\
Department of Higher Mathematics\\
Moscow Institute of Physics and Technology\\
141700 Dolgoprudny\\ RUSSIA\\
Department of Composition\\
Moscow Tchaikovsky Conservatory\\
125009 Moscow\\
RUSSIA}

\maketitle

\begin {abstract}
Using the residue class rings we analyze the structure of
sequences of musical intervals for the all-interval tone-semitone
series and the quart modes as well. The relations we stated show
the enharmonic return of the series to the initial point.
\end {abstract}

\section {Introduction}

It seems that for the first time the musical intervals were
considered as mathematical transformations of the residue class
ring modulo 12 in the paper by french scientist Camille Durutte in
1855. Notice that the modular arithmetics (the residue class
rings) was introduced by the famous German mathematicien Carl
Friedrich Gauss in 1801. Such the point of view acquired  a
special urgency in the XX century because the appearance of the
serial techniques and the set theory especially brightly
developed in the works by Milton Babbit (\cite {Andr}). It is
connected with the fruitful applications of such an approach to
constructing the sequencies of sounds by means of the certain
group of transformations of the ring ${\mathbb Z}/12{\mathbb Z}$
possessing some mathematical symmetry.

In our sketch we shall consider from this point of view the
all-interval tone-semitone series (ATsS). The scientific analysis
of this object was given by one of the authors of the present
paper (\cite {Gol2, Gol3}) and Irina Severina (\cite {Sev}). We
shall compare ATsS with the quart (non-octave) modes including the
ancient Russian "custom mode", his modern prolonged version
introduced by Yuri Butsko in Polyphonian Concert (1969) and the
artificial "synagogue mode" (really, the custom mode being
transformed) introduced by Alfred Shnitke in his 4th symphony
(1984). Notice that the phenomenon of ATsS was of interest for
the other researchers in different contexts (\cite {Voi1, Voi2,
Kul, Shar}).

\section {The all-interval tone-semitone series}

Let us enumerate a sequence of sounds of the imaginary "infinite"
keyboard by integer numbers $n=0,\pm 1,\pm 2,\dots ,$ such that
each two neighboring numbers are corresponding to the sounds
differing by a semitone. The musical interval of a duration $m$
can be considered as the operation transforming the set of
integer numbers by the rule $n\to n+m$. Because the feature of a
human ear the tones corresponding to the numbers $n$ and $n+12$
will be perceived as sounding in unison. Adding or subtracting
from each a number $n$ the appropriate integer number divisible
by 12 it is possible to transpose all the sounds of the keyboard
to one octave corresponding to the integer numbers $n=0,1,2,\dots
,11$.

Denote by ${\mathbb Z}/12{\mathbb Z}=\{0,1,2,\dots ,11\}$ the set
of sounds of this octave. If we take into account a value of
interval up to adding an integer number of octaves, then the
interval is a set of 11 transformations of the elements of the
set ${\mathbb Z}/12{\mathbb Z}$ determined by a mathematical
operation of the "adding modulo 12". This operation can be defined
as follows. The result of the adding of the numbers $n$ and $m$
modulo 12 is the number $n+m+12k$ (the notation is $m+n\ mod\
12$), where the integer $k$ is picked up such that $m+n+12k\in
{\mathbb Z}/12{\mathbb Z}$. Analogously one can define the
"multiplication modulo 12". The result of this operation is the
number $mn+12k$, where $k$ is chosen such that $mn+12k\in {\mathbb
Z}/12{\mathbb Z}$. The set ${\mathbb Z}/12{\mathbb Z}$ with the
operations of the adding modulo 12 and the multiplication modulo
12 is said to be the ring of residue classes modulo 12.

Thus, to obtain a certain row of sounds (in which we shall take
into account not only what sounds are included in it but the
order of their appearance in the row also) one should define an
operation $T$ on the ring ${\mathbb Z}/12{\mathbb Z}$. In the
author's papers \cite {Gol1, Gol2} ATsS was defined (with a
possibility to continue to the "infinite" keyboard) such that the
n-th sound of the row differs from the first on $\frac
{n(n-1)}{2}$ semitones, that is the 12th sound differs from the
first on 66 semitones. Remember that ATsS consists of 12 sounds
because under "all intervals" of its basic form occupying five
and half octaves (66:12) we mean the intervals from 1 to 11
contained inside the octave and having a musical structural
sense. In other cases, it is possible to talk about all the
intervals from 1 to 6 (because the intervals 6-11 are their
conversions inside the octave), from 0 to 24 or the imaginary
"infinite" all-interval row (within the limits of audibility and
besides it). It should be especially mentioned that the intervals
from 0 to 12 (24,36,...) are unisons (repetitions and octave
doublings of the sounds) and they have nor musical structural but
only mathematical sense. Continuing the row we get that the
intervals from 13 to 23 gives us a tritone transposition of the
basic form of ATsS which identically coincides with its
retrograde inversion, the intervals from 25 to 35 is a repetition
of sounds in the basic form and so on.

Denote $|n>$ the element of the ring ${\mathbb Z}/12{\mathbb Z}$
corresponding to the $(n-1)$th sound of ATsS. Then, the operation
$V$ determining the series by the formula $|n>=Vn,\ n=0,1,2,\dots$
can be defined by means of the following recurring formula,
\begin {equation}\label {F1}
V0=0,\ Vn=|n-1>+n\ mod\ 12,\ n=1,2,3,\dots
\end {equation}
The solution to (\ref {F1}) is given by the formula of the
arithmetic progression:
\begin {equation}\label {F2}
|n>=\frac {n(n+1)}{2}\ mod\ 12,\ n=0,1,2,\dots
\end {equation}
Notice that it follows from (\ref {F2}) that
\begin {equation}\label {F3}
|n+12>=|n>+6\ mod\ 12.
\end {equation}
At first, it implies that $|n+24k>=|n>,\ k=1,2,\dots $, i.e. the
basic form of ATsS generates a series of 24 elements of the ring
${\mathbb Z}/12{\mathbb Z}$ being repeated cyclically with the
period equal to 24. On the other hand, the operation $T$
determined on the ring ${\mathbb Z}/12{\mathbb Z}$ by the formula
$$
Tn=n+6\ mod\ 12
$$
is a tritone transposition from the musical point of view. Thus,
(\ref {F3}) implies that all the sounds from 13 to 24 obtained by
a subsequent adding to the basic form of ATsS the intervals
increasing by means of the same principle are tritone
transpositions (coinciding with retrograde inversions) of the
first 12 sounds. Notice that the continued in two times
("doubled") form of ATsS returns enharmonically to the initial
point after 17.5 octaves (210:12) of the imaginary "infinite"
keyboard and consists of 24 sounds (12x2).

Then, the operation $R_{m}$ defined on the row $x$ consisting of
$m$ elements $[x_{1},\dots ,x_{m}]$ of the ring ${\mathbb
Z}/12{\mathbb Z}$ by the formula
$$
R_{m}x_{s}=x_{m-s+1}
$$
results in a retrograde inversion of the initial row. It follows
from (\ref {F2}) that
\begin {equation}\label {F4}
R_{6}|n>=|n>+6n+3\ mod\ 12.
\end {equation}
On the other hands, the same formula implies that
\begin {equation}\label {F5}
|n+6>=|n>+6n+9\ mod\ 12.
\end {equation}
Thus, formulae (\ref {F4}) and (\ref {F5}) give us the relation
$$
|n+6>=R_{6}|n>+6\ mod\ 12
$$
confirming that in the basic form of ATsS the sounds $7-12$ are
obtained by means of the tritone transposition of the retrograde
inversion of the sounds $1-6$.

Notice that the ring ${\mathbb Z}/12{\mathbb Z}$ includes zero
divisors, i.e. the non-zero elements $m$ and $n$ such that $mn=0\
mod\ 12$. Namely, there are two such the elements, that are $3$
and $4$ because
$$
3\cdot 4=12=0\ mod\ 12.
$$
The existence of the zero divisors determines a reiteration of
sounds in the series of $12$ elements of the basic form of ATsS.
In fact, the identity
$$
|n>=|m>
$$
is equivalent to
$$
(n-m)(n+m+1)=24k
$$
for a certain integer number $k$. Given integer numbers $m$ and
$n$, one of the numbers $n-m$ or $n+m+1$ is even and the other is
odd. If $n-m$ is even, then we get the product of two integer
numbers of the form $\frac {n-m}{2}(n+m+1)=12k$, and the zero
divisors in ${\mathbb Z}/12{\mathbb Z}$ are $\frac {n-m}{2}$ and
$n+m+1$. Analogously, if $n+m+1$ is even, then $(n-m)\frac
{n+m+1}{2}=12k$ is a product of two integer numbers, and the zero
divisors are $n-m$ and $\frac {n+m+1}{2}$. Substituting as zero
divisors the pair $(3,4)$ and the pair $(1,12)$ as well, we
obtain that $12$ sounds of the basic form of ATsS contain $4$
pairs of reiterated sounds (remark that although $12$ is not a
zero divisor but it is zero in the ring ${\mathbb Z}/12{\mathbb
Z}$).

So, $12$ sounds of the basic form of ATsS contain $4$ pairs of
reiterated sounds and $4$ non-reiterated sounds, that are $8$
different sound pitches at all. For the whole chromatic scale one
needs $4$ sounds more. Why? Let us consider the equation
\begin {equation}\label {F6}
|n>=m\ mod\ 12.
\end {equation}
The equation (\ref {F6}) in the ring ${\mathbb Z}/12{\mathbb Z}$
is equivalent to the condition
\begin {equation}\label {F7}
n^{2}+n=2m+24k.
\end {equation}
The number on the right hand side of the equality (\ref {F7}) is
always even but the number in the left hand side is even only if
$n$ is even or $n=1$ or $n=11$. Hence, the equation (\ref {F6}) is
solvable only for $8$ different values of $m$.

To obtain failing $4$ sounds one can use the inversion $I$ acting
to the elements of the ring ${\mathbb Z}/12{\mathbb Z}$ by the
formula
$$
In=12-n\ mod\ 12.
$$
Following down from the first sound of the basic form of ATsS let
us construct $12$ sounds of its inversion form by the formula
$$
|-n>=I|n>,\ n=0,1,\dots ,11.
$$
Then, the equation $|-n>=m\ mod\ 12$ is equivalent to the
condition
\begin {equation}\label {F8}
n^{2}+n=2(12-m)+24k.
\end {equation}

The formula (\ref {F8}) as well as the formula (\ref {F7}) is
solvable only in $8$ cases but the solutions to (\ref {F7}) are
$m=0,1,3,4,6,7,9,10$ while the solutions to (\ref {F8}) are
$Im=0,2,3,5,6,8,9,11$ which include the failing numbers
$m=2,5,8,11$.

\section {The quart modes}

Now let us proceed to the mathematical formulae of the ancient
Russian "custom mode" and the artificial "synagogue mode" being
constructed by quarts which are filled by seconds (the synagogue
mode was introduced by Alfred Shnitke in his $4$th symphony
(\cite {Gol1, Hol})). In particular, it will be interesting for us
to compare them with the formulae of ATsS obtained above. The
interval structure of the custom mode is the following
$$
(5=)2+2+1,\ 2+2+1,\ 2+2+1,\dots (large),\ or
$$
$$
(5=)2+1+2,\ 2+1+2,\ 2+1+2,\dots (small),\ or
$$
$$
(5=)1+2+2,\ 1+2+2,\ 1+2+2,\dots (reduced),
$$
and for the "synagogue mode" is
$$
(5=)1+1+3,\ 1+1+3,\ 1+1+3,\dots ,\ or
$$
$$
(5=)1+3+1,\ 1+3+1,\ 1+3+1,\dots ,\ or
$$
$$
(5=)3+1+1,\ 3+1+1,\ 3+1+1,\dots
$$

Denote $|n>$ the series of the elements of the ring ${\mathbb
Z}/12{\mathbb Z}$ generated by the interval rows of these modes.
In the both cases, we get
\begin {equation}\label {F9}
|n+3>=|n>+5\ mod\ 12.
\end {equation}
The least common multiple of the numbers $5$ and $12$ is
$60=5\cdot 12$. In this way, it follows from (\ref {F9}) that the
continued quart modes enharmonically return to "the initial
point" through $5$ octaves $(60:12)$ and consist of $36$ sounds
$(3\cdot 12)$. The formula (\ref {F9}) shows a mathematical
symmetry in the construction of the interval rows considered as
well as the formula (\ref {F3}) determines a regularity of the
enharmonic return to "the initial point"  on the imaginary
"infinite" keyboard of the doubled form of ATsS.

To investigate the quart modes it is useful to involve the ring
${\mathbb Z}/5{\mathbb Z}=\{0,1,2,3,4\}$ in which, analogously to
the case of the ring ${\mathbb Z}/12{\mathbb Z}$, the operations
of the "adding modulo $5$" and the "multiplication modulo $5$" are
introduced. Notice that the number $5$ (unlike $12$) is prime.
Hence the ring ${\mathbb Z}/5{\mathbb Z}$ doesn't contain zero
divisors. It follows that ${\mathbb Z}/5{\mathbb Z}$ is a field
(in spite of the operations of adding, subtraction and
multiplication the operation of division is defined).

Let us denote the series of elements of the field ${\mathbb
Z}/5{\mathbb Z}$ generated by the interval rows of the continued
(large, small and reduced) custom and sinagogue modes in each of
the cases considered above by $|Ln>,|Sn>,|Rn>$ and $|SGn>$,
respectively. In all the cases we obtain three elements of the
field ${\mathbb Z}/5{\mathbb Z}$, namely,
$$
|L0>=0,\ |L1>=2,\ |L2>=4,
$$
$$
|S0>=0,\ |S1>=2,\ |S2>=3,
$$
$$
|R0>0,\ |R1>=1,\ |R2>=3,
$$
$$
|SG0>=0,\ |SG1>=1,\ |SG2>=2>,\ or
$$
$$
|SG0>=0,\ |SG1>=1,\ |SG2>=4,\ or
$$
$$
|SG0>=0,\ |SG1>=3,\ |SG2>=4.
$$
It is easy to see that the sequences considered consisting of
triples of the elements of the field ${\mathbb Z}/5{\mathbb Z}$
use up all possible constructions which can be generated by the
quart modes determined by the formula (\ref {F9}). Notice that
from the viewpoint of mathematics the large custom mode possesses
the most common symmetry because the elements of the field
${\mathbb Z}/5{\mathbb Z}$ which are generated by it form a
cyclic group with the generator $2$:
$$
|Ln>=2n\ mod\ 5,\ n=0,1,2.
$$

\section {Discussions}

We have derived the basic mathematical formulae of the
all-interval tone-semitone series and the quart modes including
the ancient Russian custom mode and the artificial "synagogue
mode" proposed by Alfred Shnitke as well. The information
contained in these formulae allows to know when the row will
enharmonically return to the "initial point" and what sounds does
it consist of also. The musical sense of some mathematical
operations on the residue rings is revealed.

\end {document}